\documentclass[preprint]{elsarticle}       % Specifies the document class
\usepackage{lipsum}
\makeatletter
\def\ps@pprintTitle{%
 \let\@oddhead\@empty
 \let\@evenhead\@empty
 \def\@oddfoot{}%
 \let\@evenfoot\@oddfoot}
\makeatother

\usepackage{lineno,hyperref}
\usepackage{rotating}
\usepackage{makeidx}
\usepackage{float}
\usepackage{epsfig}
\usepackage{amsmath,amssymb}
\usepackage[latin1]{inputenc}
\usepackage{epic,eepic}
\usepackage{overpic}
\usepackage{color}
\usepackage{amsfonts}
\usepackage{graphicx}
\usepackage{mwe} % For dummy images
\usepackage{subcaption}

\usepackage{algorithm}
\usepackage{algorithmicx}
\usepackage{algpseudocode}
\usepackage{enumitem}

\newenvironment{demo}{\smallskip\noindent{{\it Proof.}}\hskip \labelsep}%
            {\hfill\penalty10000\raisebox{-.09em}{\large\bf\rm $\blacksquare$}\par\medskip}

\newtheorem{theorem}{Theorem}[section]

\newtheorem{proposition}[theorem]{Proposition}
\newtheorem{corollary}[theorem]{Corollary}
\newtheorem{definition}{Definition}
\newtheorem{example}{Example}[section]

\newcommand{\cF}{\mathcal{F}}
\DeclareMathOperator*{\argmin}{arg\,min}

\def\qno{\mathcal{Q}^{\text{NL}}}

\journal{}

\begin{document}

\small

\begin{frontmatter}
%Weighted Essentially Non-Oscillatory 1D-Moving Least Squares based on Partition of Unity method
\title{Integrating Moving Least Squares with non-linear WENO method:\\ A novel Partition of Unity approach in 1D}\tnotetext[label1]{Dionisio F. Y\'a\~nez has been supported through project CIAICO/2021/227 (Proyecto financiado por la Conselleria de Innovaci\'on, Universidades, Ciencia y Sociedad digital de la Generalitat Valenciana), by grant PID2020-117211GB-I00 and by PID2023-146836NB-I00 funded by MCIN/AEI/10.13039/501100011033.}

\author[UV]{Inmaculada Garc\'es}
\ead{Inmaculada.Garces@uv.es}
\author[UV]{Jos\'e M. Ram\'on}
\ead{Jose.Manuel.Ramon@uv.es}
\author[UPCT]{Juan Ruiz-\'Alvarez}
\ead{juan.ruiz@upct.es}
\author[UV]{Dionisio F. Y\'a\~nez}
\ead{Dionisio.Yanez@uv.es}

\date{Received: date / Accepted: date}

\address[UV]{Departamento de Matem\'aticas. Universidad de Valencia, Valencia (Spain).}
\address[UPCT]{Departamento de Matem\'atica Aplicada y Estad\'istica. Universidad  Polit\'ecnica de Cartagena, Cartagena (Spain).}

%The correct dates will be entered by the editor.

\begin{abstract}
%The approximation of noisy data is a classical problem that appears in several applications as computer-aided geometric design, numerical
%resolution of partial differential equations, and design of curves and surfaces. Some methods have been developed in this way with interesting
%results when the data does not present any discontinuity. Moving least squares is an adequate strategy to fit some inputs and it has been employed
%satisfactorily in statistics and applied mathematics. However, when the data presents an isolated discontinuity some non-desired effects occur as Gibbs phenomenon.
%In this paper, we mix moving least squares with the well-known non-linear technique WENO in order to construct a non-linear operator that improves the
%approximations close to the discontinuities and conserves the order of accuracy in the smooth zones. We study its properties in several dimensions. Some numerical experiments are presented to check the theoretical results.

The approximation of data is a fundamental challenge encountered in various fields, including computer-aided geometric design, the numerical solution of partial differential equations, or the design of curves and surfaces. Numerous methods have been developed to address this issue, providing good results when the data is continuous. Among these, the Moving Least Squares (MLS) method has proven to be an effective strategy for fitting data, finding applications in both statistics and applied mathematics. However, the presence of isolated discontinuities in the data can lead to undesirable artifacts, such as the Gibbs phenomenon, which adversely affects the quality of the approximation.

In this paper, we propose a novel approach that integrates the Moving Least Squares method with the well-established non-linear Weighted Essentially Non-Oscillatory (WENO) method. This combination aims to construct a non-linear operator that enhances the accuracy of approximations near discontinuities while maintaining the order of accuracy in smooth regions. We investigate the properties of this operator in one dimension, demonstrating its effectiveness through a series of numerical experiments that validate our theoretical findings.
\end{abstract}
\end{frontmatter}
%All acknowledgements should be placed in the back of the paper after Conclusions..

\section{Introduction}

In the context of numerical approximation methods, we introduced in \cite{bsplines} a novel WENO B-spline-based quasi-interpolation algorithm. This algorithm relayed on a non-linear modification of the B-spline basis functions, which form a partition of unity. This innovative approach not only preserved the smoothness of the original spline but also adapted effectively to discontinuities in the underlying function, reducing Gibbs-like oscillations close to those discontinuities. The application of WENO weights to the B-spline functions is an original idea that can be applied to other constructions where a partition of unity is part of the quasi-interpolation operator. Starting from this point, the current paper presents a natural continuation of our earlier research by exploring a similar non-linear modification within the framework of the MLS method for univariate data. In this work, as in our previous article \cite{bsplines}, we exploit the partition of unity in the quasi-interpolation operator to make it nonlinear. This new approach integrates the MLS method with the non-linear WENO technique, enhancing the accuracy near discontinuities while maintaining the performance in smooth regions.%Just as we applied a partition of unity in our previous work, this new approach integrates the MLS method with the non-linear WENO technique, further enhancing the accuracy of approximations near discontinuities while maintaining performance in smooth regions.

The MLS is a general algorithm designed to approximate a functional from some scattered data (see \cite{davidlevin}). It is an important technique used in several contexts with different denominations. Thus, in statistics, MLS is called {\it local polynomial regression}  \cite{FASSHAUER,LOADER}. This is an interesting procedure because it can be used in applications such as denoising, image processing, subdivision schemes, and to solve numerically partial differential equations (see for example \cite{carretal,davidlevin2}). Following the notation presented in \cite{davidlevin}, we consider a normed function space $\cF$ on $\mathbb{R}^n$ and let $f\in\cF$ and $\{L_i(f)\}_{i=1}^N$ be some scattered data points, where $\{L_i\}_{i=1}^N$ are bounded linear functionals on $\cF$. The problem consists in approximating $L(f)$ from the known data, where $L$ is another functional. Typically, this value is computed as a linear combination of the data, i.e.
\begin{equation}\label{eq1}
\mathcal{Q}(f)=\sum_{i=1}^N a_i L_i(f),
\end{equation}
where $\{a_i\}_{i=1}^N$ are values calculated imposing the following condition: Let $P\equiv \text{span}\{p_j\}_{j=1}^J\subseteq \mathcal{F}$ be a finite set of fundamental functions (usually polynomials) then
\begin{equation}\label{eq2}
\mathcal{Q}(p)=\sum_{i=1}^N a_i L_i(p)=L(p), \quad p\in P.
\end{equation}
From \eqref{eq1} and \eqref{eq2} we can analyze the error. Let $\Omega_0\subset \mathbb{R}^n$ be the support of the functional $L$, $\Omega_I\subset \mathbb{R}^n$ be the support of $\sum_{i=1}^na_iL_i$, and we define $\Omega=\Omega_0\cup\Omega_I$ then:
\begin{equation}
\begin{split}
|L(f)-\mathcal{Q}(f)|&= |L(f)-L(p)+L(p)-\mathcal{Q}(f)|=|L(f)-L(p)+\mathcal{Q}(p)-\mathcal{Q}(f)|\\
&\leq (\|L\|+\|\mathcal{Q}\|)\|f-p\|_\Omega\leq (\|L\|+\sum_{i=1}^N |a_i|\|{L_i}\|)\|f-p\|_\Omega.
\end{split}
\end{equation}
In \cite{davidlevin}, we find the following definition that we reproduce here in order to make this paper almost self-contained.
\begin{definition}[The MLS approximation]
Let $\{L_i\}_{i=1}^N$, $L$ and $P$ defined as above, and let $\{\Theta(L_i,L)\}$ be some non-negative weights. The MLS approximation to $L(f)$ is defined as $L(p^*)$, where $p^*\in P$ is minimizing, among all $p\in P$, the weighted least-squares error functional
\begin{equation}\label{problema1}
\sum_{i=1}^N(L_i(p)-L_i(f))^2\Theta(L_i,L).
\end{equation}
\end{definition}
However, in many real-world applications, data can exhibit discontinuities or strong gradients, making a linear treatment inappropriate for approximating \( L(f) \). This is due to Gibbs-like oscillations and smearing near these discontinuities. The WENO method has been developed and analysed in recent years, and there is a substantial amount of literature on this topic (see, e.g., \cite{JiangShu, Liu, doi:10.1137/070679065}). WENO employs a non-linear combination of multiple Lagrange interpolations of a specific degree to achieve a higher order of accuracy, provided that the data is continuous and free of discontinuities. The key feature of the WENO method is its ability to identify when the stencils used for interpolation encounter a non-smooth region in the data, relying on specific operators known as smoothness indicators.

In this paper, we integrate the WENO and MLS methods to create a non-linear operator that effectively avoids the Gibbs phenomenon and provides accurate approximations when the data is smooth. For this purpose, we partition our data set $\{L_i(f)\}_{i=1}^N$ into different subsets $\Upsilon_k$, such that
$$
\{L_i(f)\}_{i=1}^N = \bigcup_k \Upsilon_k,
$$
and define the non-linear operator:
$$
\mathcal{Q}^\text{NL}(f) = \sum_{k} \alpha_k(\Upsilon_k) \mathcal{Q}^k(f),
$$
where $\mathcal{Q}^k(f)$ is the solution to the problem \eqref{problema1} for the data $\Upsilon_k$, and $\alpha_k(\Upsilon_k)$ are non-linear weights defined using the values of $\Upsilon_k$ for each $k$.

Throughout the paper, we consider the point values data discretization, where the functionals $L_i$ represent evaluations of a function at scattered data points $\{x_i\}_{i=1}^N$. The results can be easily extended to the cell-average discretization \cite{harten3}, commonly used in image and signal processing.
%In this context, we introduce a kind of data discretization that will be referenced throughout the paper, we consider the approximation of point values, where the functionals $L_i$ represent the evaluations of a function at scattered data points $\{x_i\}_{i=1}^N$. In a similar way, we can extend the result to the cell-average discretization, which is commonly used in the context of image and signal processing.

\subsection{MLS for point-value (PV) data}\label{mlspv}
We start focusing our work on the classical least squares approximation example. We assume that $\{x_i\}_{i=1}^N$ are some points in $[a,b]\subseteq \mathbb{R}$ and $L_i(f)=f(x_i)$, %+\epsilon_i$
$i=1,\hdots,N$ is a sampling of the function $f$ at the sites $x_i$. %where $\epsilon_i$ is some noise (it can be zero),
Then, we try to approximate $L_x(f)=f(x)$ with $x\in[a,b]$, minimizing the following problem in the set of polynomials of degree less than or equal to $m$, $\Pi_m$:
\begin{equation}\label{problema2}
p_x=\argmin_{p\in\Pi_m}\sum_{i=1}^N (p(x_i)-L_i(f))^2w(|x-x_i|).
\end{equation}
Here, $w$ is a non-negative weight function with specific properties that we analyze later in Subsection \ref{PVmls}. %$d_\text{PV}$ denotes the distance between the  point $\bx$, where we want to approximate, and the point $\bx_i$, both of which are in $\mathbb{R}^n$.
Now, we can define the quasi-interpolation operator $\mathcal{Q}(f)(x) = p_x(x)$, where $\mathcal{Q}(f)$ is given by,
%where $\theta$ is a non-negative weight function with some properties that we will analyze in Section \ref{}, $d_\text{PV}$ is a distance between two points in $\mathbb{R}^n$, and evaluate it in $\bx$, i.e., $\mathcal{Q}(f)(\bx)=p_\bx(\bx)$. It is well known that the resulting operator $\mathcal{Q}$ is independent of the data $L_i(f)$ and linear. In fact, we can represent
\begin{equation}\label{q}
\mathcal{Q}(f)(x)=\sum_{i=1}^N a_i(x)L_i(f),
\end{equation}
where $\{a_i\}_{i=1}^N$ are real functions with the same degree of smoothness as the weight function $w$. It is well known that the resulting operator $\mathcal{Q}$ is independent of the data $L_i(f)$ and linear.  Also, due to the definition of the problem in \eqref{problema2}, the linear operator $\mathcal{Q}$ reproduces polynomials in $\Pi_m$, i.e.
$$\mathcal{Q}(p)=p,\quad \forall\,p\in\Pi_m.$$
%These two properties are very important to obtain an approximation to data obtained sampling a continuous function.
%When discontinuous data is approximated by a continuous and regular function the Gibbs phenomenon appears close to the discontinuities \cite{gottshu}.
%To construct a non-linear operator capable of avoiding Gibbs phenomenon and to approximate adequately when the data belongs to a smooth part, we follow the method introduced in the previous section: We divide our domain, $\Omega$ in subdomains $\Upsilon_k$, with
%$\Omega=\cup_k \Upsilon_k$ and define the non-linear operator:
%$$\qno(f)(\bx)=\sum_{k}\alpha_k(\Upsilon_k) p^k_\bx(\bx),$$
%being $p_\bx^k$ the solution of the problem \eqref{problema2} in $\Upsilon_k$, and $\alpha_k(\Upsilon_k)$ some non-linear weights defined depending on the values of
%$\Upsilon_k$.
These two properties are very important for obtaining an accurate approximation of a smooth function through data in the point values. On the other hand, it is well known that when data that present discontinuities is approximated using a continuous and regular function, the Gibbs phenomenon appears close to the discontinuities \cite{gottshu}. To construct a non-linear operator capable of avoiding the Gibbs phenomenon and, at the same time, adequately approximating the data at smooth parts, we follow the method introduced in the previous section: We divide our domain, $[a,b]$, into subdomains $[a_k,b_k]$, such that
$$
[a,b] \subseteq \bigcup_k [a_k,b_k],\quad x\in[a_k,b_k]\,\,\forall \, k,
$$
and define the non-linear operator:
$$
\qno(f)(x) = \sum_{k} \alpha_k([a_k,b_k]) p^k_x(x),
$$
where $p_x^k$ is the solution to the problem in \eqref{problema2} in $[a_k,b_k]$, and $\alpha_k([a_k,b_k])$ are non-linear weights defined based on the values of $[a_k,b_k]$. We will discuss in Section \ref{WENO_MLS} how to design these non-linear weights.

\subsection{WENO method}
The WENO algorithm has emerged as a powerful tool for the approximation of data values, particularly in scenarios where discontinuities are present (see e.g. \cite{doi:10.1137/070679065}). Originally developed for solving hyperbolic partial differential equations, the WENO method is designed to maintain high accuracy while preventing spurious osci\-llations that can arise near discontinuities \cite{JiangShu}. This is achieved through the use of adaptive weights that emphasize smoother regions of the data while minimizing the influence of stencils that cross the discontinuities.

Let us consider a function $f$ that we wish to approximate at the point $x$. The WENO reconstruction at this point is formulated as:

\begin{equation}
\hat{f}(x) = \sum_{k=0}^{r-1} \beta_k p_k(x),
\end{equation}
where $r$ indicates the number of stencils, $p_k$ are the polynomial approximations derived from each stencil, and $\beta_k$ represents the non-linear weights. These weights are determined by:

\begin{equation}
\beta_k = \frac{\alpha_k}{\sum_{l=0}^{r-1} \alpha_l},
\,\,\text{with}\,\,
\alpha_k = \frac{C_k}{(\epsilon + I_k)^t},
\end{equation}
where $C_k$ denotes the linear weights, $I_k$ are the smoothness indicators, $\epsilon$ is a small positive constant to prevent division by zero, and $t$ is typically set to 2 to optimize accuracy in smooth regions. The smoothness indicators $I_k$ are designed to quantify the smoothness of the data, often inspired by the concept of total variation:

\begin{equation}
I_k = \sum_{l=1}^{r-1} \int_{x_{i-1/2}}^{x_{i+1/2}} \Delta x^{2l-1} \left( \frac{d^l}{dx^l} f_k(x) \right)^2 dx,
\end{equation}
where $\Delta x$ represents the grid spacing. Considering the specific challenges of our problem, we will utilize these concepts to develop a non-linear quasi-interpolation operator for the point values discretization, while also redefining the smoothness indicators in Section \ref{WENO_MLS} to better suit our context.

\subsection{Outline of this paper}

The rest of the paper is divided as follows: in Section \ref{PVmls} we review the solution for the linear MLS problem and present the explicit form for the point values setting in 1D. The new WENO-like MLS method will be introduced in Section \ref{WENO_MLS} and some theoretical results will be proved. We specially focus our attention on the Gibbs phenomenon, which appears close to the discontinuities in the linear case and it is avoided in the non-linear one, and on the order of the approximation in the smooth parts of the data. Some experiments are performed in section \ref{NE} in order to check the theoretical results. Finally, in Section \ref{conc}, we provide some conclusions and remarks.

\section{Reviewing the general solution of MLS problem}\label{PVmls}

%In this section, we employ the results obtained in \cite{davidlevin} to solve the MLS problem introduced in \eqref{problema1}. Thus, we suppose that $P= \text{span}\{p_j\}_{j=1}^J$, also that $\{L_i(f)\}_{i=1}^N$ with $1\in P$ and $J\leq N$. Now we construct the matrices $E$ and $D$ as
%\begin{equation}\label{construccionE}
%E_{i,j}=L_i(p_j), \,\, 1\leq i \leq N, 1\leq j\leq J, \quad D=\text{diag}(\Theta(L_1,L),\hdots,\Theta(L_N,L)).
%\end{equation}
%We assume that $\text{rank}(E)=J$, and get the solution given by:
%\begin{equation}\label{solucionMLS}
%\mathcal{Q}(f)=(L_1(f),\hdots,L_N(f))DE(E^TDE)^{-1}(L(p_1),\hdots,L(p_J))^T.
%\end{equation}
%This general solution can be applied to the two kinds of data discretizations introduced in Subsections \ref{mlspv} and \ref{mlsca}. The next two subsections are dedicated to this purpose.

In this section, we employ the results obtained in \cite{davidlevin} to solve the MLS problem introduced in \eqref{problema1}. We assume that $P = \text{span}\{p_j\}_{j=1}^J$ and that $\{L_i(f)\}_{i=1}^N$ includes $1 \in P$ with $J \leq N$. We construct the matrices $E$ and $D$ as follows:
\begin{equation}\label{construccionE}
E_{i,j} = L_i(p_j), \quad 1 \leq i \leq N, \quad 1 \leq j \leq J, \quad D = \text{diag}(\Theta(L_1,L), \ldots, \Theta(L_N,L)).
\end{equation}
We assume that $\text{rank}(E) = J$, leading to the solution given by:
\begin{equation}\label{solucionMLS}
\mathcal{Q}(f) = (L_1(f), \ldots, L_N(f)) D E (E^T D E)^{-1} (L(p_1), \ldots, L(p_J))^T.
\end{equation}
This general solution can be applied to some kinds of data discretizations. In particular, in the next subsection we apply it to the point-value case.

\subsection{Point-value MLS in 1d}

Let $\{x_i\}_{i=1}^N$ be a set of points in the interval $[a,b]$ and $f$ an unknown function. We use the point value discretization, so that $L_i(f)=f(x_i)$, $i=1,\hdots,N$. We consider a value $x\in[a,b]$ and we want to approximate $L(f)=f(x)$ using polynomials up to degree $d$. Thus, we employ the basis
$$p_j(t)=(t-x)^j, \quad j=0,\hdots,d.$$
Then $E_{i,j}=L_i(p_j)=(x_i-x)^j$ with rank$(E)=d+1$ since we have a Vandermonde matrix. To design $\Theta(L_i,L)$, we define a weight function using a function $w:[0,\infty)\to[0,1]$ with the following properties:
\begin{enumerate}
\item We choose $w$ to be a decreasing function, with $\lim_{x\to\infty} w(x)=0.$
\item The function at the point $x=0$ is 1, i.e. $w(0)=1$.
\item There exists $p\in\mathbb{N}$ such that $w\in \mathcal{C}^p([0,\infty))$.
\end{enumerate}
Therefore we define $\Theta(L_i,L):=w_i(x):=w(|x-x_i|/h)$, where $h:=\max_{i=2,\hdots,N}\{|x_{i}-x_{i-1}|\}$ is the fill distance (see e.g. \cite{FASSHAUER}). We consider some classical examples described in \cite{FASSHAUER} and collected in Table \ref{tabla1nucleos}, 
 where the cutoff function $(\cdot)_+:\mathbb{R} \to \mathbb{R}$ is defined as:
 \begin{equation*}
 (x)_+=\begin{cases}
 x, & x\geq 0,\\
 0, & x<0.
 \end{cases}
 \end{equation*}

 \begin{table}[h!]
\centering
\begin{tabular}{lll}
\hline
$w(  \,r)$ & RBF &                                                          \\ \hline
$e^{-  r^2}$ & Gaussian $\mathcal{C}^\infty$ & G                                          \\
$\left(1 + r^2 \right)^{-1/2}$ & Inverse MultiQuadratic $C^\infty$ & IMQ  \\
$e^{-  r}$ & Mat\'ern $\mathcal{C}^0$ & M0 \\
$e^{-  r} \left( 1 +   r \right)$ & Mat\'ern $\mathcal{C}^2$ & M2 \\
$e^{-  r} \left( 3 + 3  r +   r^2 \right)$ & Mat\'ern $\mathcal{C}^4$ & M4 \\
$(1 -   r)^2_+$ & Wendland $\mathcal{C}^0$ & W0 \\
$(1 -   r)^4_+ \left( 4  r + 1 \right)$ & Wendland $\mathcal{C}^2$ & W2 \\
$(1 -   r)^6_+ \left( 35  r^2 + 18  r + 3 \right)$ & Wendland $\mathcal{C}^4$ & W4 \\ \hline
\end{tabular}
\caption{Examples of RBFs.}\label{tabla1nucleos}
\end{table}
In \cite{davidlevin}, the author applies the MLS with $w(x)=e^{-x^2}$. In this last case, if {$w_i(x)<\varepsilon$, then we suppose that $w_i(x)=0$}, being $\varepsilon$ a small parameter, in our experiments $\varepsilon=10^{-9}$. Also other functions introduced in \cite{LOADER} or \cite{lopezyanez} with the form:
$$w(x)=\begin{cases}
(1-x^q)^p, & x\in [0,1],\\
0, & \text{otherwise,}
\end{cases}
$$
with $p,q\in\mathbb{N}$ can be used. With these ingredients, we apply the formula described in \eqref{solucionMLS} and as $p_j(x)=0$ for all $j\neq 0$, we get:
\begin{equation}\label{solMLSpv}
\mathcal{Q}(f)(x)=\sum_{i=1}^N(DE(E^TDE)^{-1})_{i,1}f(x_i).
\end{equation}
Now we can denote $C_{i}(x)=(DE(E^TDE)^{-1})_{i,1}$. To provide an explicit form of this formula and prove the following basic propositions, we will follow the same method demonstrated in \cite{lopezyanez}.%and give an explicit form of this formula, proving the following basic propositions, following the same way showed in \cite{lopezyanez}.
\begin{proposition}
With the same notation as before, $C_{i_0}$ has the explicit form
$$C_{i_0}(x)=\frac{w_{i_0}(x)}{|E^TDE|}\sum_{i_1,\hdots,i_d=1}^N w_{i_1}(x)\hdots w_{i_d}(x)(x_{i_1}-x)(x_{i_2}-x)^2\hdots (x_{i_d}-x)^d\prod_{0\leq k<l \leq d}(x_{i_l}-x_{i_k}), \quad 1\leq i_0 \leq N.$$
\end{proposition}
\begin{demo}
By simplicity, we write $w_i(x)$ as $w_i$. From
$$
E^TD\mathbf{f}=
\begin{bmatrix}
w_1 &  w_2 & \hdots & w_N\\
w_1(x_1-x) & w_2(x_2-x) & \hdots & w_N(x_N-x)\\
\vdots & \vdots& \ddots & \vdots\\
w_1(x_1-x)^d & w_2(x_2-x)^d & \hdots & w_N((x_N-x)^d
\end{bmatrix}\begin{bmatrix}
f_1\\
f_2\\
\vdots\\
f_N
\end{bmatrix}=
\begin{bmatrix}
\sum_{i=1}^Nw_if_i\\
\sum_{i=1}^Nw_if_i(x_i-x)\\
\vdots\\
\sum_{i=1}^Nw_if_i(x_i-x)^d
\end{bmatrix},
$$
and from the expression
 $$E^TDE=\begin{bmatrix}
\sum_{i=1}^Nw_i &  \sum_{i=1}^N w_i(x_i-x) & \hdots & \sum_{i=1}^N w_i(x_i-x)^d\\
\sum_{i=1}^Nw_i(x_i-x) &   \sum_{i=1}^N w_i(x_i-x)^2 &\hdots & \sum_{i=1}^N w_i(x_i-x)^{d+1}\\
\vdots & \vdots& \ddots & \vdots\\
\sum_{i=1}^Nw_i(x_i-x)^d & \sum_{i=1}^Nw_i(x_i-x)^{d+1} & \hdots & \sum_{i=1}^N w_i(x_i-x)^{d+d}
\end{bmatrix},
$$
we have that
\begin{equation*}
\begin{split}
|E^TDE|&=
\begin{vmatrix}
\sum_{i_0=1}^Nw_{i_0} &  \sum_{i_1=1}^N w_{i_1}(x_{i_1}-x) & \hdots & \sum_{i_d=1}^N w_{i_d}(x_{i_d}-x)^d\\
\sum_{i_0=1}^Nw_{i_0}(x_{i_0}-x) &   \sum_{i_1=1}^N w_{i_1}(x_{i_1}-x)^2 &\hdots & \sum_{i_d=1}^N w_{i_d}(x_{i_d}-x)^{d+1}\\
\vdots & \vdots& \ddots & \vdots\\
\sum_{i_0=1}^Nw_{i_0}(x_{i_0}-x)^d & \sum_{i_1=1}^Nw_{i_1}(x_{i_1}-x)^{d+1} & \hdots & \sum_{i_d=1}^N w_{i_d}(x_{i_d}-x)^{d+d}
\end{vmatrix}\\
&=\sum_{i_0=1}^N w_{i_0}\begin{vmatrix}
1 &  \sum_{i_1=1}^N w_{i_1}(x_{i_1}-x) & \hdots & \sum_{i_d=1}^N w_{i_d}(x_{i_d}-x)^d\\
(x_{i_0}-x) &   \sum_{i_1=1}^N w_{i_1}(x_{i_1}-x)^2 &\hdots & \sum_{i_d=1}^N w_{i_d}(x_{i_d}-x)^{d+1}\\
\vdots & \vdots& \ddots & \vdots\\
(x_{i_0}-x)^d & \sum_{i_1=1}^Nw_{i_1}(x_{i_1}-x)^{d+1} & \hdots & \sum_{i_d=1}^N w_{i_d}(x_{i_d}-x)^{d+d}
\end{vmatrix}\\
&=\sum_{i_0,i_1,\hdots,i_d=1}^N w_{i_0}w_{i_1}\hdots w_{i_d}\begin{vmatrix}
1 &  (x_{i_1}-x) & \hdots & (x_{i_d}-x)^d\\
(x_{i_0}-x) &   (x_{i_1}-x)^2 &\hdots & (x_{i_d}-x)^{d+1}\\
\vdots & \vdots& \ddots & \vdots\\
(x_{i_0}-x)^d & (x_{i_1}-x)^{d+1} & \hdots & (x_{i_d}-x)^{d+d}
\end{vmatrix}\\
&=\sum_{i_0,i_1,\hdots,i_d=1}^N w_{i_0}w_{i_1}\hdots w_{i_d}(x_{i_1}-x)(x_{i_2}-x)^2\hdots (x_{i_d}-x)^d\begin{vmatrix}
1 &  1 & \hdots & 1\\
(x_{i_0}-x) &   (x_{i_1}-x) &\hdots & (x_{i_d}-x)\\
\vdots & \vdots& \ddots & \vdots\\
(x_{i_0}-x)^d & (x_{i_1}-x)^{d} & \hdots & (x_{i_d}-x)^{d}
\end{vmatrix}\\
&=\sum_{i_0,i_1,\hdots,i_d=1}^N w_{i_0}w_{i_1}\hdots w_{i_d}(x_{i_1}-x)(x_{i_2}-x)^2\hdots (x_{i_d}-x)^d\prod_{0\leq k<l \leq d}(x_{i_l}-x_{i_k}).\\
%&=\sum_{i_0=1}^N w_{i_0}\sum_{i_1,\hdots,i_d=1}^N w_{i_1}\hdots w_{i_d}(x_{i_1}-x)(x_{i_2}-x)^2\hdots (x_{i_d}-x)^d\prod_{0\leq k<l \leq d}(x_{i_l}-x_{i_k})\\
\end{split}.
\end{equation*}
By Kramer's formula, we have that
\begin{equation*}
\begin{split}
&\begin{vmatrix}
\sum_{i_0=1}^nw_{i_0}f_{i_0} &  \sum_{i_1=1}^n w_{i_1}(x_{i_1}-x) & \hdots & \sum_{i_d=1}^n w_{i_d}(x_{i_d}-x)^d\\
\sum_{i_0=1}^nw_{i_0}(x_{i_0}-x)f_{i_0} &   \sum_{i_1=1}^n w_{i_1}(x_{i_1}-x)^2 &\hdots & \sum_{i_d=1}^n w_{i_d}(x_{i_d}-x)^{d+1}\\
\vdots & \vdots& \ddots & \vdots\\
\sum_{i_0=1}^nw_{i_0}(x_{i_0}-x)^df_{i_0} & \sum_{i_1=1}^nw_{i_1}(x_{i_1}-x)^{d+1} & \hdots & \sum_{i_d=1}^n w_{i_d}(x_{i_d}-x)^{d+d}
\end{vmatrix}=\\
&=\sum_{i_0=1}^n w_{i_0}f_{i_0}\begin{vmatrix}
1 &  \sum_{i_1=1}^n w_{i_1}(x_{i_1}-x) & \hdots & \sum_{i_d=1}^n w_{i_d}(x_{i_d}-x)^d\\
(x_{i_0}-x) &   \sum_{i_1=1}^n w_{i_1}(x_{i_1}-x)^2 &\hdots & \sum_{i_d=1}^n w_{i_d}(x_{i_d}-x)^{d+1}\\
\vdots & \vdots& \ddots & \vdots\\
(x_{i_0}-x)^d & \sum_{i_1=1}^nw_{i_1}(x_{i_1}-x)^{d+1} & \hdots & \sum_{i_d=1}^n w_{i_d}(x_{i_d}-x)^{d+d}
\end{vmatrix}\\
&=\sum_{i_0,i_1,\hdots,i_d=1}^n f_{i_0}w_{i_0}w_{i_1}\hdots w_{i_d}(x_{i_1}-x)(x_{i_2}-x)^2\hdots (x_{i_d}-x)^d\prod_{0\leq k<l \leq d}(x_{i_l}-x_{i_k})\\
&=\sum_{i_0=1}^n w_{i_0}f_{i_0}\sum_{i_1,\hdots,i_d=1}^n w_{i_1}\hdots w_{i_d}(x_{i_1}-x)(x_{i_2}-x)^2\hdots (x_{i_d}-x)^d\prod_{0\leq k<l \leq d}(x_{i_l}-x_{i_k}).\\
\end{split}.
\end{equation*}
Therefore,
$$C_{i_0}(x)=\frac{n_{i_0}(x)}{|E^TDE|}=\frac{w_{i_0}}{|E^TDE|}\sum_{i_1,\hdots,i_d=1}^N w_{i_1}\hdots w_{i_d}(x_{i_1}-x)(x_{i_2}-x)^2\hdots (x_{i_d}-x)^d\prod_{0\leq k<l \leq d}(x_{i_l}-x_{i_k}).$$
\end{demo}
We show a basic example to clarify the formulas. We can see that the smoothness of the operator only depends on the smoothness of the function $w(|\cdot|)$.
\begin{example}
 Let be $x_i=i-1$, $i\in\mathbb{Z}$, we will take $n$ points and construct
the approximation for any $x\in[0,n-1]$ (included the boundary values) for polynomials of degree $d=1$ and $n=3$, then,
$$E=
\begin{bmatrix}
1 & -x  \\
1 & 1-x \\
1 & 2-x \\
\end{bmatrix},
\quad
D=
\begin{bmatrix}
w_1(x) & 0&  0\\
0 & w_2(x) &  0\\
0 & 0 &  w_3(x)
\end{bmatrix},
$$

\begin{equation*}
\begin{split}
|E^TDE|
=&\sum_{i_0,i_1=1}^3 w_{i_0}w_{i_1}(x_{i_1}-x)\prod_{0\leq k<l \leq 1}(x_{i_l}-x_{i_k})=\sum_{i_0,i_1=1}^3 w_{i_0}w_{i_1}(x_{i_1}-x)(x_{i_1}-x_{i_0})\\
=&w_{1}\sum_{i_1=1}^3 w_{i_1}(x_{i_1}-x)(x_{i_1}-x_{1})+w_{2}\sum_{i_1=1}^3 w_{i_1}(x_{i_1}-x)(x_{i_1}-x_{2})+w_{3}\sum_{i_1=1}^3 w_{i_1}(x_{i_1}-x)(x_{i_1}-x_{3})\\
=&w_{1}w_{2}(x_{2}-x)(x_{2}-x_{1})+w_{1}w_{3}(x_{3}-x)(x_{3}-x_{1})+w_{2}w_{1}(x_{1}-x)(x_{1}-x_{2})+w_{2}w_{3}(x_{3}-x)(x_{3}-x_{2})\\
&+w_{3}w_{1}(x_{1}-x)(x_{1}-x_{3})+w_{3}w_{2}(x_{2}-x)(x_{2}-x_{3})\\
%=&w_{1}w_{2}(x_{2}-x)+2w_{1}w_{3}(x_{3}-x)+w_{2}w_{1}(x-x_{1})+w_{2}w_{3}(x_{3}-x)+2w_{3}w_{1}(x-x_{1})+w_{3}w_{2}(x-x_{2})\\
%=&w_{1}w_{2}(x_{2}-x+x-x_{1})+2w_{1}w_{3}(x_{3}-x+x-x_{1})+w_{2}w_{3}(x_{3}-x+x-x_{2})\\
=&w_{1}w_{2}+4w_{1}w_{3}+w_{2}w_{3}.\\
\end{split}
\end{equation*}
Analogously,
\begin{equation*}
\begin{split}
n_1(x)=&w_{1}\sum_{i_1=2}^3 w_{i_1}(x_{i_1}-x)(x_{i_1}-x_{1})=w_{1}w_{2}(x_{2}-x)(x_{2}-x_{1})+w_{1}w_{3}(x_{3}-x)(x_{3}-x_{1})\\
=&w_{1}w_{2}(x_{2}-x)+2w_{1}w_{3}(x_3-x)=w_{1}w_{2}+4w_{1}w_{3}+(-2w_{1}w_{3}-w_{1}w_{2})x,\\
n_2(x)=&w_{2}\sum_{i_1=1,i_1\neq 2}^3 w_{i_1}(x_{i_1}-x)(x_{i_1}-x_{2})=w_{2}w_{1}(x_{1}-x)(x_{1}-x_{2})+w_{2}w_{3}(x_{3}-x)(x_{3}-x_{2})\\
=&w_{2}w_{1}x+w_{2}w_{3}(x_{3}-x)=2w_{2}w_{3}+(w_{2}w_{1}-w_{2}w_{3})x,\\
n_3(x)=&w_{3}\sum_{i_1=2}^3 w_{i_1}(x_{i_1}-x)(x_{i_1}-x_{3})=w_{3}w_{1}(x_{1}-x)(x_{1}-x_{3})+w_{3}w_{2}(x_{2}-x)(x_{2}-x_{3})\\
=&-w_3w_2+(2w_3w_1+w_3w_2)x.\\
\end{split}
\end{equation*}
Therefore,
\begin{equation*}
\begin{split}
C_1(x)&=\frac{w_{1}w_{2}+4w_{1}w_{3}+(-2w_{1}w_{3}-w_{1}w_{2})x}{w_{1}w_{2}+4w_{1}w_{3}+w_{2}w_{3}},\\
C_2(x)&=\frac{2w_{2}w_{3}+(w_{2}w_{1}-w_{2}w_{3})x}{w_{1}w_{2}+4w_{1}w_{3}+w_{2}w_{3}},\\
C_3(x)&=\frac{-w_3w_2+(2w_3w_1+w_3w_2)x}{w_{1}w_{2}+4w_{1}w_{3}+w_{2}w_{3}}.
\end{split}
\end{equation*}
\end{example}

Note that, in this work, we could use the concept of $h-\rho-\delta$ sets of mesh size $h$, density $\leq \rho$  and separation $\delta$ (see \cite{davidlevin}) but it is sufficient to consider that there exists a constant $C$ such that:
$$h=\max_{i=2,\hdots,N}\{|x_{i}-x_{i-1}|\} \leq Ch_m=C\min_{i=2,\hdots,N}\{|x_{i}-x_{i-1}|\},$$
i.e., using the notation by Wendland in \cite{wendland} (see definition 4.6 of \cite{wendland}) the set $\{x_i\}_{i=1}^N$ of data sites is quasi-uniform.

Finally, we review the following theorem proved in \cite{davidlevin}, and also mentioned in \cite{FASSHAUER}, where the order of approximation is proved when the data comes from the discretization of a sufficiently continuous function.

\begin{theorem}\label{repropolis}
Let $(a,b) \subset \mathbb{R}$ be an open set, $x\in (a,b)$, $\{x_i\}_{i=1}^N \subset (a,b)$ a set of $N$ distinct nodes and  $\{f_i=f(x_i)\}_{i=1}^N$ a corresponding set of function values with $f\in\Pi_d(\mathbb{R}^n)$. Then the MLS approximation defined in Eq. \eqref{solMLSpv} satisfies
$$\mathcal{Q}(f)(x)=f(x).$$
\end{theorem}
From this Theorem, Th. \ref{repropolis}, we get the order of accuracy.
\begin{corollary}\label{cororder}
Let $(a,b) \subset \mathbb{R}$. If $f\in\mathcal{C}^{d+1}(a,b)$, $\{x_i\}_{i=1}^N \subset (a,b)$ are quasi-uniformly distributed nodes with fill distance $h$, the weight function $w$ is compactly supported with support size $c$, then the approximation defined in Eq. \eqref{solMLSpv}  fulfills
$$|\mathcal{Q}(f)(x)-f(x)|\leq Ch^{d+1}\max_{\xi\in [a,b]}|f^{(d+1)}(\xi)|, \quad x\in \Omega,$$
where $C$ is a constant independent of $h$.
\end{corollary}

%Analogously, it is possible to prove the same results for cell-average setting. 
Now, we have the sufficient ingredients to construct the non-linear operator based on the partition of unity method (see \cite{FASSHAUER}).

\section{WENO-like MLS based on the partition of unity method}\label{WENO_MLS}
In this section, we explain the non-linear method based on partition of unity. The idea is quite simple, to divide the domain into several subdomains, to measure the ``smoothness'' of the data in each subdomain, and to use only those which are free of discontinuities.

 We suppose an open domain $\Omega \subset \mathbb{R}$ with $[a,b]\subset \Omega$,  some data points $\{x_i\}_{i=1}^N$ with $h$ the fill distance, and some associated set $\{f_i=f(x_i)\}_{i=1}^N$ with $f:\Omega \to \mathbb{R}$ an unknown function. Let $x\in [a,b]$
be a point of the domain, and we suppose that the compact support of $w$ is $c$. We consider that $\{\tilde{x}_k\}_{k=1}^m\subseteq [a,b]$ are some center points and call
\begin{equation}\label{omegak}
\Omega_k=\{x\in[a,b]: \tilde{w}_k(x)=w(\gamma_k|x-\tilde{x}_k|/h)>0\},\quad k=1,\hdots,m
\end{equation}
%\begin{equation}\label{chik}
%\chi_k(x)=\chi(x,\gamma_k)=\{i\in\mathbb{N}: w_i(x)=w(\gamma_k|x-x_i+c_kh|)>0\}=\{i\in\mathbb{N}:|x-x_i+c_kh|\leq \frac{c}{\gamma_k}h\},\quad k=1,\hdots,m
%\end{equation}
where $m$ is the number of subsets chosen to do the partition of unity, $\gamma_k>0$ is a constant to determine the number of data of each problem. Now, we solve a MLS problem for each $k$
$$p_k=\underset{p\in\Pi_d}{\arg \min} \sum_{x_i\in \Omega_k}(p(x_i)-f(x_i))^2w(\gamma_k|x-x_i|/h),$$
Note that the constant $\gamma_k$ and the centers $\{\tilde{x}_k\}_{k=1}^m$ have to be selected to satisfy that:
$$|\Omega_k\cap\{x_i\}_{i=1}^N|> d+1,\quad \forall k=1,\hdots,m,\quad \text{and}\quad \Omega\subseteq \bigcup_{k=1}^m\Omega_k$$
typically, (see e.g. \cite{cavo}) this constant is called the {\it shape parameter}.
Thus, a partition of unity method can be designed
\begin{equation}\label{PUM}
\mathcal{Q}_\text{PU}(f)(x)=\sum_{k=1}^m \theta_k (x)\mathcal{Q}_k(f)(x)=\sum_{k=1}^m \theta_k (x)p_k(x),
\end{equation}
being
$$\theta_k(x)=\frac{\delta_k(x)}{\sum_{k=1}^m\delta_j(x)}\quad \text{with} \quad \delta_k(x)=w(\gamma_k|x-\tilde{x}_k|/h).$$
Now, we solve the following LS problem:
$$\tilde{p}_k=\underset{p\in\Pi_d}{\arg \min} \sum_{x_i\in \Omega_k}(p(x_i)-f(x_i))^2.$$
It is clear by Cor. \ref{cororder} that, if the function $f\in\mathcal{C}^{d+1}$, then:
$$||\tilde{p}_k-f||_{\infty,\Omega_k}=O(h^{d+1})  \quad \text{as}\,h\to 0, \quad\forall \, {k}=1,\hdots,m.$$
However, if a discontinuity crosses one subset, $k_0$ then
$$||\tilde{p}_{k_0}-f||_{\infty,\Omega_{k_0}}=O(1) \quad \text{as}\,h\to 0.$$
Therefore, we can define the smoothness indicators as:
$$\mathcal{I}_k=\frac{1}{|\Omega_k\cap \{x_i\}_{i=1}^N|}\sum_{x_i\in \Omega_k\cap \{x_i\}_{i=1}^N}|\tilde{p}_k(x_i)-f(x_i)|,$$
where $|\Omega_k\cap \{x_i\}_{i=1}^N|$ is the number of points belonging to this set. {Thus, the smoothness indicator is defined as the mean of the absolute errors from fitting a polynomial using the least squares method. This measure helps to assess the quality of the polynomial approximation. For example, if we think about gridded data with grid size $h$ and a polynomial of degree $n$, the mean error would be $O(h^{n+1})$ when the data is smooth. If there's a jump discontinuity, the mean error would be $O(1)$, indicating difficulty in capturing the discontinuity. This difference in error behavior helps to detect and quantify data smoothness.}

With these ingredients we can write the non-linear 1D- MLS operator:
\begin{equation}\label{NLPUMLS}
\mathcal{Q}^\text{NL}_\text{PU}(f)(x)=\sum_{k=1}^m \beta_k (x)\mathcal{Q}_k(f)(x)=\sum_{k=1}^m \beta_k (x)p_k(x),
\end{equation}
with
\begin{equation}\label{NLPUMLS2}
\beta_k(x)=\frac{\alpha_k(x)}{\sum_{j=1}^m\alpha_j(x)}, \quad \alpha_k(x)=\frac{\theta_k(x)}{\mathcal{I}^t_k+\epsilon }, \quad k=1,\hdots,m,
\end{equation}
with $\epsilon$ a constant to avoid zero in the denominator and $t$ a parameter to assure the optimal accuracy in the smooth regions. We use $t=4$ for all the experiments. We take in our experimental examples, $\epsilon=10^{-14}$. By Eqs. \eqref{NLPUMLS} and  \eqref{NLPUMLS2}
is clear that the smoothness of the operator $\mathcal{Q}^\text{NL}_\text{PU}(f)$ is the same as the smoothness presented by the function $w$.
\begin{theorem}
Let $[a,b] \subset \Omega \subset \mathbb{R}$ be an open set, $\{x_i\}_{i=1}^N \subset [a,b]$ are quasi-uniformly distributed with fill distance $h$, the weight function $w$ is compactly supported with support size $c$. Let $x\in [a,b]$ be a point and $\Omega_k$ the subsets defined in Eq. \eqref{omegak}, then if  $t\geq1$ and
%$$\exists\,\, k_0  \in \{1,\hdots,m\} : \sum_{i\in \chi_{k_0}(x)}|f(x_i)-p_{k_0}(x_i)|=O(h^{d+1})$$
 $$\exists\,\, k_0  \in \{1,\hdots,m\} : \mathcal{I}_{k_0}=O(h^{d+1})$$
 then
 $$|\mathcal{Q}^\text{NL}_\text{PU}(f)(x)-f(x)|=O(h^{d+1}).$$
with $\mathcal{Q}^\text{NL}_\text{PU}(f)$ defined in Eq. \eqref{NLPUMLS}.
\end{theorem}
\begin{demo}
We know that if $t\geq1$ then
$$\alpha_k(x)=\begin{cases}
O(h^{-(d+1)}), & \text{if} \,\, f  \,\, \text{is smooth in} \,\, [x_{\min\chi_k(x)},x_{\max\chi_k(x)}],\\
O(1), & \text{if} \,\, f  \,\, \text{is non-smooth in} \,\, [x_{\min\chi_k(x)},x_{\max\chi_k(x)}].\\
\end{cases}
$$
Thus, $\sum_{k\in \chi_k(x)} \alpha_k=O(h^{-(d+1)})$ and
$$\beta_k(x)=\begin{cases}
O(1), & \text{if} \,\, f  \,\, \text{is smooth in} \,\, [x_{\min\chi_k(x)},x_{\max\chi_k(x)}],\\
O(h^{d+1}), & \text{if} \,\, f  \,\, \text{is non-smooth in} \,\, [x_{\min\chi_k(x)},x_{\max\chi_k(x)}].\\
\end{cases}
$$
 We denote as $K^x_1=\{k\in \mathbb{N}: f  \,\, \text{is smooth in} \,\, [x_{\min\chi_k(x)},x_{\max\chi_k(x)}]\}$ and we get
\begin{equation*}
\begin{split}
|\mathcal{Q}^\text{NL}_\text{PU}(f)(x)-f(x)|&=\left|\sum_{k\in \chi_k(x)}\beta_k(x) p_k(x)-f(x)\right|=\left|\sum_{k\in \chi_k(x)}\beta_k(x) p_k(x)-\sum_{k\in \chi_k(x)}\beta_k(x) f(x)\right|\\
&=\left|\sum_{k\in \chi_k(x)}\beta_k(x) (p_k(x)- f(x))\right|\\
&\leq
\sum_{k\in \chi_k(x)}\beta_k(x) \left|p_k(x)- f(x)\right|\\
&=\sum_{k\in K^x_1}\beta_k(x) \left|p_k(x)- f(x)\right|+\sum_{k\notin K^x_1}\beta_k(x) \left|p_k(x)- f(x)\right|\\
&=\sum_{k\in K^x_1}O(1)O(h^{d+1})+\sum_{k\notin K^x_1}O(h^{d+1})O(1)\\
&=O(h^{d+1}).
\end{split}
\end{equation*}
\end{demo}

\section{Numerical experiments}\label{NE}
%%Todos los programas están en la carpeta:Articulos_dioni/wenolikeMLS/programas/ultima_version

%In this section, we study the different characteristics of the new algorithm and check the theoretical properties. For this, we design some experiments divided in subsections, we start with an example  to analyze the order of accuracy using an smooth function. After that, we use the new method to approximate non-continuous function and verify that the non desired phenomenons close to the discontinuities are avoided.

In this section, we examine the various characteristics of the new algorithm and verify its theoretical properties. To do this, we design several experiments, divided into two subsections. We begin with an example to analyze the order of accuracy using a smooth function. Next, we apply the new method to approximate a non-continuous function and confirm that undesired phenomena near discontinuities are avoided.

\subsection{Order of accuracy}
We start this section analyzing the order of accuracy. For this reason, we start with a smooth function
$$f(x)=\sin(\pi x),\,\, x\in \mathbb{R},$$
discretized in the interval $[-3,3]$ with an uniform grid $N=2^l+1$, $h_l=2^{-l}$, $\{x^l_i=-3+\frac{3}{2^{l-1}}i\}_{i=0}^{2^l}$,  our data are $\{f^l_i=f(x_i^l)\}_{i=0}^{2^l}$. We calculate an approximation at the points $z_j=\frac{j}{1000}$, $j=0,\hdots,1000$, and get the numerical errors
\(e^l_j=|f(z_j)-\mathcal{I}^l(z_j)|\), where $\mathcal{I}^l$ is $\mathcal{Q}_\text{PU}(f)$ and $\mathcal{Q}^\text{NL}_\text{PU}(f)$ in each level $l$.
\begin{equation}\label{erroryorden}
\begin{split}
\text{MAE}_l=\max_{j=0,\hdots,1000}e^l_j,  \quad r_l^{\infty}=\frac{\log(\text{MAE}_{l-1}/\text{MAE}_{l})}{\log(h_{l-1}/h_{l})},
\end{split}
\end{equation}
%we will use the acronyms MLSPU$^p_{\mathcal{H}}$ and NL-MLSPU$^p_{\mathcal{H}}$ when the operator shown in Eqs. \eqref{PUM} and \eqref{NLPUMLS} are used with $p$ the degree of polynomials and $\mathcal{H}=\text{G,W2,W4}$ when $w$ is Gaussian, Wendland $\mathcal{C}^2$ or Wendland $\mathcal{C}^4$ functions showed in Table \ref{tabla1nucleos}. 
we will use the acronyms MLSPU$^p_{\mathcal{H}}$ and NL-MLSPU$^p_{\mathcal{H}}$ when the operators shown in Eqs. \eqref{PUM} and \eqref{NLPUMLS} are used, with $p$ representing the degree of the polynomials and $\mathcal{H} = \text{G, W2, W4}$ indicating that $w$ is a Gaussian, Wendland $\mathcal{C}^2$, or Wendland $\mathcal{C}^4$ function, as shown in Table \ref{tabla1nucleos}. In this subsection, for all $k$, the shape parameter $\gamma_k=0.15$ for the W2 and W4, and $\gamma_k=0.7$ for G. Finally, the center points chosen are the data points, i.e.,
$\{\tilde{x}^l_k\}_{k=0}^{2^l}=\{{x}^l_i\}_{i=0}^{2^l}$. In this first example, we can see in Table \ref{exp1} that the numerical order of accuracy for $p=2,3$ is 4, the results are very similar in both cases, both when using the linear method and the nonlinear method.
\begin{table}[!ht]
\begin{center}
\begin{tabular}{lrrrrrrrrrrr}
& \multicolumn{2}{c}{MLSPU$^2_{\text{W2}}$} & &   \multicolumn{2}{c}{NL-MLSPU$^2_{\text{W2}}$}& &\multicolumn{2}{c}{MLSPU$^3_{\text{W2}}$} & &   \multicolumn{2}{c}{NL-MLSPU$^3_{\text{W2}}$}\\ \cline{1-3} \cline{5-6} \cline{8-9} \cline{11-12}
$l$ & $\text{MAE}_l$ & $r_l^{\infty}$ & &$\text{MAE}_l$ & $r_l^{\infty}$ & & $\text{MAE}_l$ & $r_l^{\infty}$ & &$\text{MAE}_l$ & $r_l^{\infty}$       \\
\hline
$7$  &                   4.0219e-04 &        & &  7.7337e-04&         &&        3.6743e-04 &        &&   1.1402e-04 &                   \\
$8$  &                   2.5460e-05 &  3.9816& &  7.7214e-05&   3.3242&&        2.3249e-05 &  3.9822&&   2.3186e-05 &  2.2980       \\
$9$  &                   1.5964e-06 &  3.9954& &  1.5964e-06&   5.5960&&        1.4576e-06 &  3.9955&&   1.4576e-06 &  3.9916       \\
$10$ &                   9.9855e-08 &  3.9988& &  9.9855e-08&   3.9988&&        9.1172e-08 &  3.9989&&   9.1172e-08 &  3.9989       \\
\hline
\hline
& \multicolumn{2}{c}{MLSPU$^2_{\text{W4}}$} & &   \multicolumn{2}{c}{NL-MLSPU$^2_{\text{W4}}$}& &\multicolumn{2}{c}{MLSPU$^3_{\text{W4}}$} & &   \multicolumn{2}{c}{NL-MLSPU$^3_{\text{W4}}$}\\ \cline{1-3} \cline{5-6} \cline{8-9} \cline{11-12}
$l$ & $\text{MAE}_l$ & $r_l^{\infty}$ & &$\text{MAE}_l$ & $r_l^{\infty}$ & & $\text{MAE}_l$ & $r_l^{\infty}$ & &$\text{MAE}_l$ & $r_l^{\infty}$       \\
\hline
$7$  &                2.6063e-04  &        & & 3.7839e-04 &        &&        2.5310e-04 &        & &  1.2709e-04 &               \\
$8$  &                1.6459e-05  & 3.9851 & & 3.3346e-05 &  3.5043&&        1.5981e-05 &  3.9852& &  1.5966e-05 &  2.9928      \\
$9$  &                1.0314e-06  & 3.9962 & & 1.0314e-06 &  5.0148&&        1.0014e-06 &  3.9962& &  1.0014e-06 &  3.9948      \\
$10$ &                6.4508e-08  & 3.9990 & & 6.4508e-08 &  3.9990&&        6.2633e-08 &  3.9990& &  6.2633e-08 &  3.9990      \\
\hline
\hline
& \multicolumn{2}{c}{MLSPU$^2_{\text{G}}$} & &   \multicolumn{2}{c}{NL-MLSPU$^2_{\text{G}}$}& &\multicolumn{2}{c}{MLSPU$^3_{\text{G}}$} & &
\multicolumn{2}{c}{NL-MLSPU$^3_{\text{G}}$}\\ \cline{1-3} \cline{5-6} \cline{8-9} \cline{11-12}
$l$ & $\text{MAE}_l$ & $r_l^{\infty}$ & &$\text{MAE}_l$ & $r_l^{\infty}$ & & $\text{MAE}_l$ & $r_l^{\infty}$ & &$\text{MAE}_l$ & $r_l^{\infty}$       \\
$7$    &           6.0703e-05 &        &&   6.0697e-05 &        &&         6.0697e-05 &        &&   5.9436e-05 &             \\
$8$    &           3.8149e-06 &  3.9920&&   3.8149e-06 &  3.9919&&         3.8146e-06 &  3.9920&&   3.8145e-06 &  3.9618    \\
$9$    &           2.3876e-07 &  3.9980&&   2.3876e-07 &  3.9980&&         2.3874e-07 &  3.9980&&   2.3874e-07 &  3.9980    \\
$10$   &           1.4928e-08 &  3.9995&&   1.4928e-08 &  3.9995&&         1.4927e-08 &  3.9995&&   1.4927e-08 &  3.9995    \\                                                                                                          \hline
\end{tabular}
\end{center}
\caption{Errors and rates using linear and non-linear MLS methods for the test function $f(x)=\sin(\pi x)$ evaluated at grid points.}\label{exp1}
\end{table}

Now, we repeat the experiment for non uniform grids using random points in the interval $[-3,3]$. We define $h_l=\max_{i=1,\hdots,2^l}\{x_i^l-x_{i-1}^l\}$ and choose  the rest of the parameters as in the previous example. We show the results in Table \ref{exp2} where we observe that the order of accuracy is the expected one,
3 for $p=2$ and 4 for $p=3$. In these cases, we can see that the results obtained for smooth data when approximating with the linear and non-linear methods are very similar. We have performed some other experiments with $p=0,1$ and the conclusions obtained are similar. Therefore, the behaviors of the new method and the linear one are analogous when working with smooth data. In the next subsection, we study the results when the algorithms are employed with data with strong gradients or discontinuities.
\begin{table}[!ht]
\begin{center}
\begin{tabular}{lrrrrrrrrrrr}
& \multicolumn{2}{c}{MLSPU$^2_{\text{W2}}$} & &   \multicolumn{2}{c}{NL-MLSPU$^2_{\text{W2}}$}& &\multicolumn{2}{c}{MLSPU$^3_{\text{W2}}$} & &   \multicolumn{2}{c}{NL-MLSPU$^3_{\text{W2}}$}\\ \cline{1-3} \cline{5-6} \cline{8-9} \cline{11-12}
$l$ & $\text{MAE}_l$ & $r_l^{\infty}$ & &$\text{MAE}_l$ & $r_l^{\infty}$ & & $\text{MAE}_l$ & $r_l^{\infty}$ & &$\text{MAE}_l$ & $r_l^{\infty}$       \\
\hline
$7$  &                   1.9596e-01 &          && 1.9235e-01  &       &&           1.7915e-01  &       &&   8.0618e-02 &                          \\
$8$  &                   4.3203e-02 &  2.4376  && 3.1909e-02  & 2.8961&&           2.8964e-02  & 2.9376&&   1.1217e-02 &  3.1795                  \\
$9$  &                   3.0308e-03 &  3.5240  && 4.0260e-03  & 2.7455&&           1.4567e-03  & 3.9653&&   2.7746e-04 &  4.9065                  \\
$10$ &                   9.2540e-04 &  2.4978  && 1.0280e-03  & 2.8743&&           2.4080e-04  & 3.7898&&   1.4554e-04 &  1.3585                  \\
$11$ &                   2.5814e-04 &  2.9947  && 1.4821e-04  & 4.5429&&           4.2180e-05  & 4.0862&&   4.1809e-05 &  2.9258                  \\
$12$ &                   2.1936e-05 &  2.9723  && 2.1712e-05  & 2.3157&&           1.7336e-06  & 3.8479&&   1.7336e-06 &  3.8373                  \\
\hline
\hline
& \multicolumn{2}{c}{MLSPU$^2_{\text{W4}}$} & &   \multicolumn{2}{c}{NL-MLSPU$^2_{\text{W4}}$}& &\multicolumn{2}{c}{MLSPU$^3_{\text{W4}}$} & &   \multicolumn{2}{c}{NL-MLSPU$^3_{\text{W4}}$}\\ \cline{1-3} \cline{5-6} \cline{8-9} \cline{11-12}
$l$ & $\text{MAE}_l$ & $r_l^{\infty}$ & &$\text{MAE}_l$ & $r_l^{\infty}$ & & $\text{MAE}_l$ & $r_l^{\infty}$ & &$\text{MAE}_l$ & $r_l^{\infty}$       \\
\hline
$7$  &               1.4197e-01 &        &&   1.4042e-01 &        &&      3.9316e-02  &       &&   3.8930e-02 &                         \\
$8$  &               3.2138e-02 &  2.3949&&   1.8806e-02 &  3.2412&&      6.6309e-03  & 2.8694&&   6.0305e-03 &  3.0065                 \\
$9$  &               1.8305e-03 &  3.8003&&   2.6980e-03 &  2.5751&&      3.0239e-04  & 4.0952&&   2.9599e-04 &  3.9977                 \\
$10$ &               7.2989e-04 &  1.9358&&   8.3613e-04 &  2.4665&&      4.8043e-05  & 3.8733&&   4.7776e-05 &  3.8399                 \\
$11$ &               2.4145e-04 &  2.5948&&   1.7455e-04 &  3.6747&&      7.5429e-06  & 4.3429&&   7.5425e-06 &  4.3300                 \\
$12$ &               1.7509e-05 &  3.1635&&   1.7446e-05 &  2.7766&&      3.3738e-07  & 3.7460&&   3.3738e-07 &  3.7459                 \\
\hline
\hline
& \multicolumn{2}{c}{MLSPU$^2_{\text{G}}$} & &   \multicolumn{2}{c}{NL-MLSPU$^2_{\text{G}}$}& &\multicolumn{2}{c}{MLSPU$^3_{\text{G}}$} & &
\multicolumn{2}{c}{NL-MLSPU$^3_{\text{G}}$}\\ \cline{1-3} \cline{5-6} \cline{8-9} \cline{11-12}
$l$ & $\text{MAE}_l$ & $r_l^{\infty}$ & &$\text{MAE}_l$ & $r_l^{\infty}$ & & $\text{MAE}_l$ & $r_l^{\infty}$ & &$\text{MAE}_l$ & $r_l^{\infty}$       \\   \hline
 $7$  &                      3.9354e-02 &         &&  3.9353e-02 &        &&       3.9294e-02 &        &&   3.8375e-02&                       \\
 $8$  &                      1.2135e-02 &  1.8967 &&  1.1948e-02 &  1.9217&&       6.6261e-03 &  2.8697&&   5.6461e-03&   3.0895              \\
 $9$  &                      1.3229e-03 &  2.9394 &&  1.3204e-03 &  2.9212&&       3.0217e-04 &  4.0952&&   2.8725e-04&   3.9501              \\
 $10$ &                      2.5344e-04 &  3.4791 &&  2.5442e-04 &  3.4671&&       4.8006e-05 &  3.8734&&   4.7557e-05&   3.7865              \\                                                                                                           $11$ &                      1.1531e-04 &  1.8471 &&  1.1438e-04 &  1.8753&&       7.5374e-06 &  4.3428&&   7.5368e-06&   4.3210              \\
 $12$ &                      8.9862e-06 &  3.0767 &&  8.9857e-06 &  3.0669&&       3.3716e-07 &  3.7459&&   3.3716e-07&   3.7458              \\
 \hline
\end{tabular}
\end{center}
\caption{Errors and rates using linear and non-linear MLS methods for the test function $f(x)=\sin(\pi x)$ evaluated at non uniform points.}\label{exp2}
\end{table}

\subsection{Avoiding oscillations close to the discontinuities}

\begin{figure}[htbp!]
\begin{center}
		\begin{tabular}{cc}
	 Linear & 	Non-linear \\
\multicolumn{2}{c}{$p=2$} \\
	\includegraphics[width=8cm, height=5.15cm]{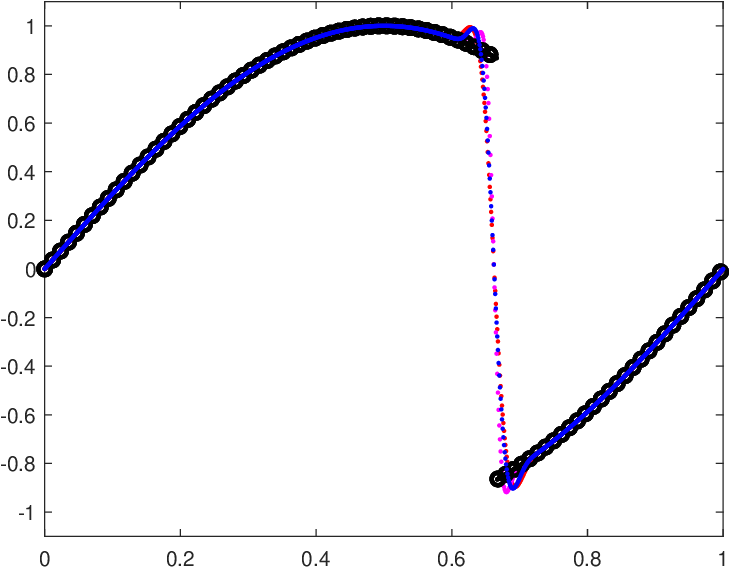} & 	\includegraphics[width=8cm, height=5.15cm]{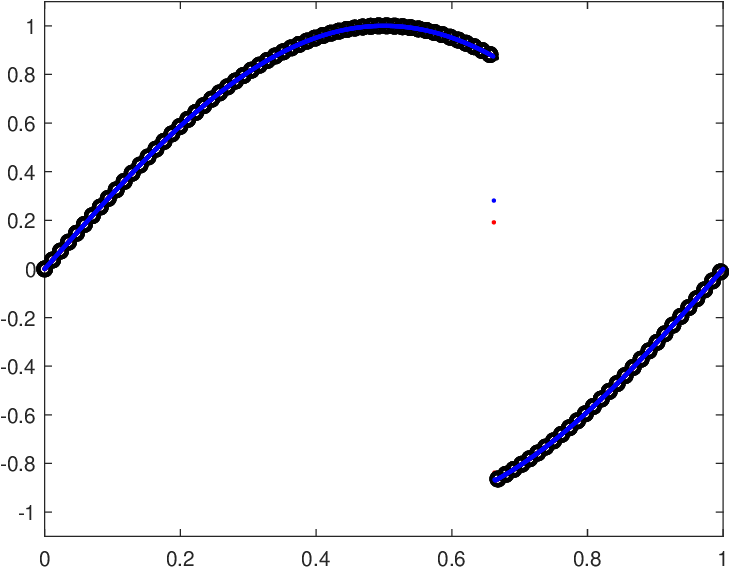}\\
	\includegraphics[width=8cm, height=5.15cm]{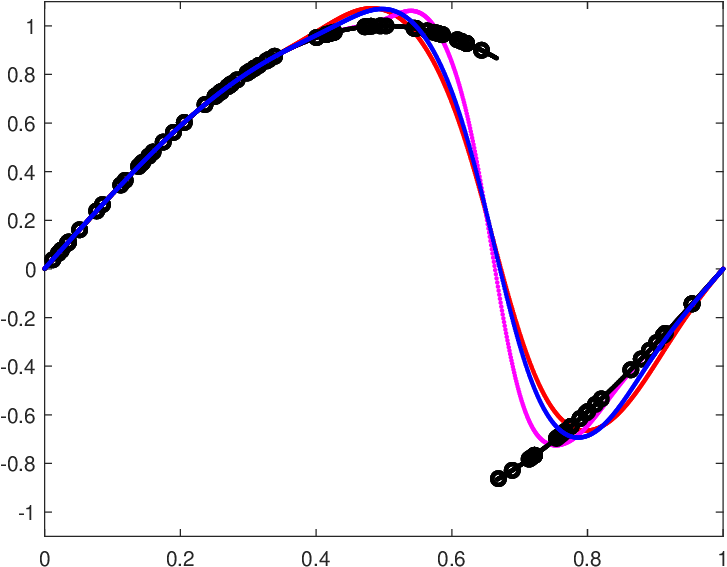} & 	\includegraphics[width=8cm, height=5.15cm]{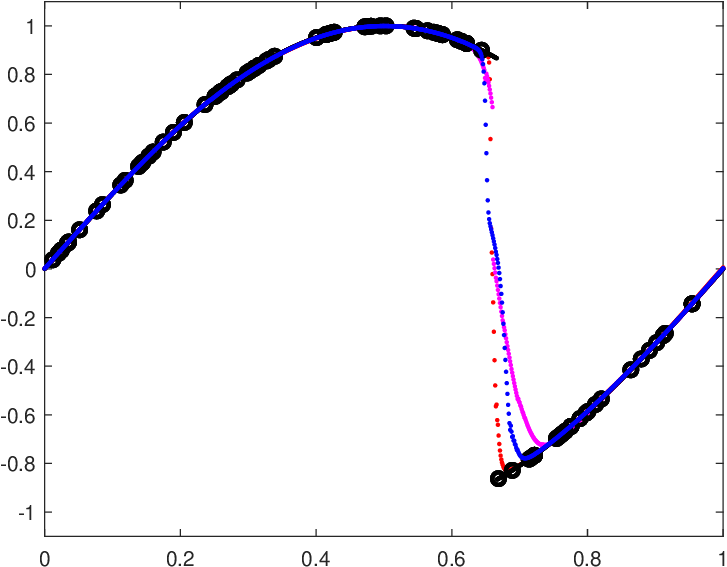}\\
\multicolumn{2}{c}{$p=3$}\\
	\includegraphics[width=8cm, height=5.15cm]{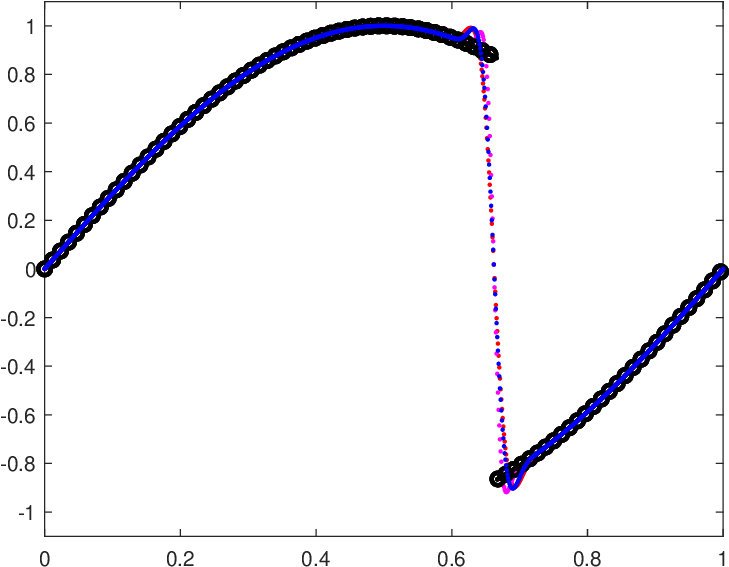} & 	\includegraphics[width=8cm, height=5.15cm]{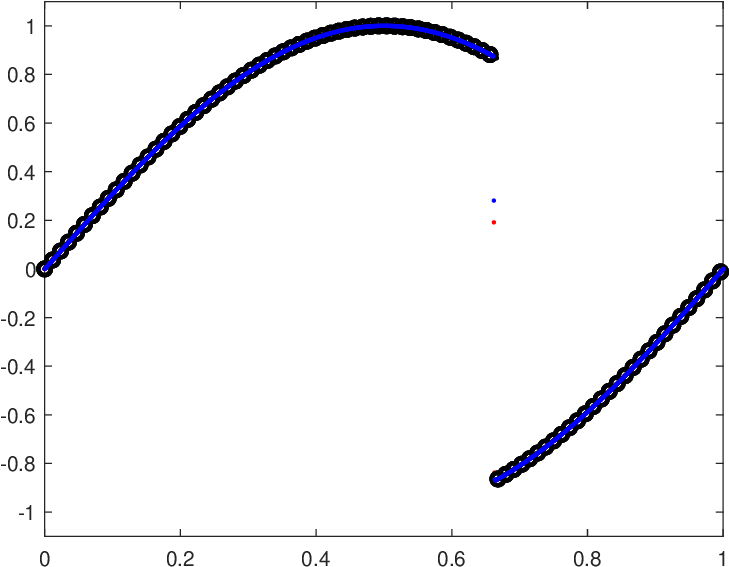}\\
	\includegraphics[width=8cm, height=5.15cm]{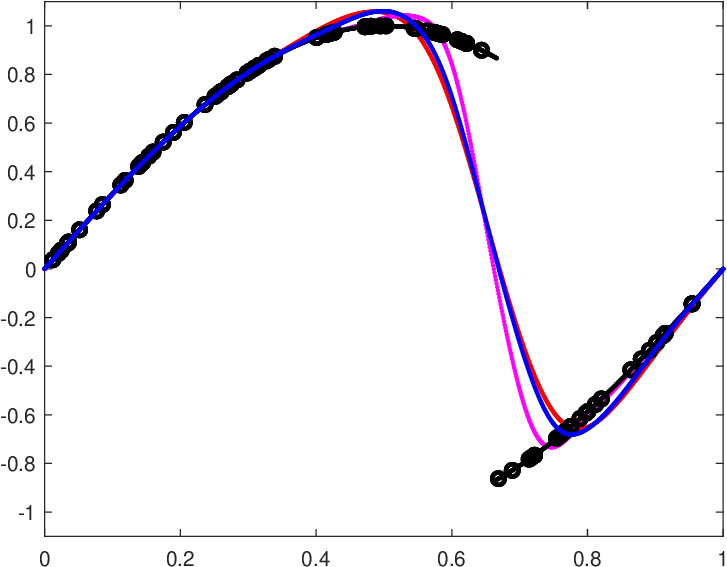} & 	\includegraphics[width=8cm, height=5.15cm]{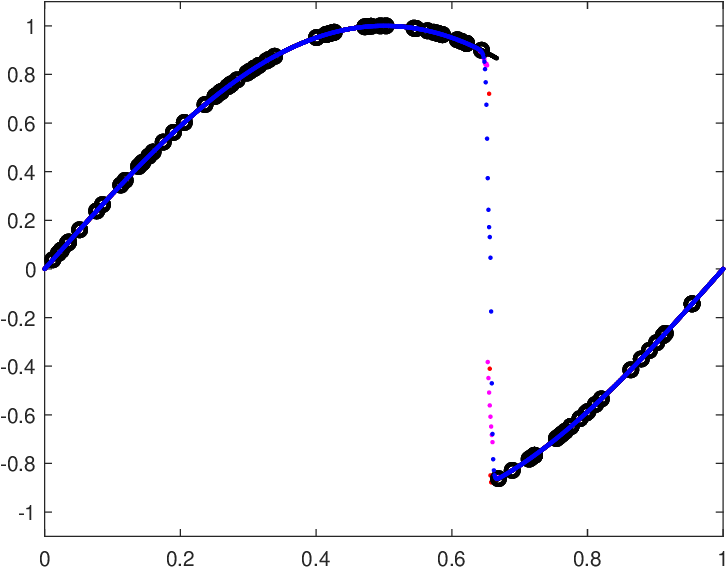}\\
		\end{tabular}
\end{center}
			\caption{Approximation to the function $g$ (black), Eq. \eqref{funciong}, using linear and non-linear methods with ${\text{W2}}$ (red), ${\text{W4}}$ (blue) and G (magenta).}
		\label{aproxfunciongw2w4}
	\end{figure}

\begin{figure}[htbp!]
\begin{center}
		\begin{tabular}{cc}
Linear & Non-linear \\
\multicolumn{2}{c}{$p=2$} \\
\includegraphics[width=8cm, height=4.15cm]{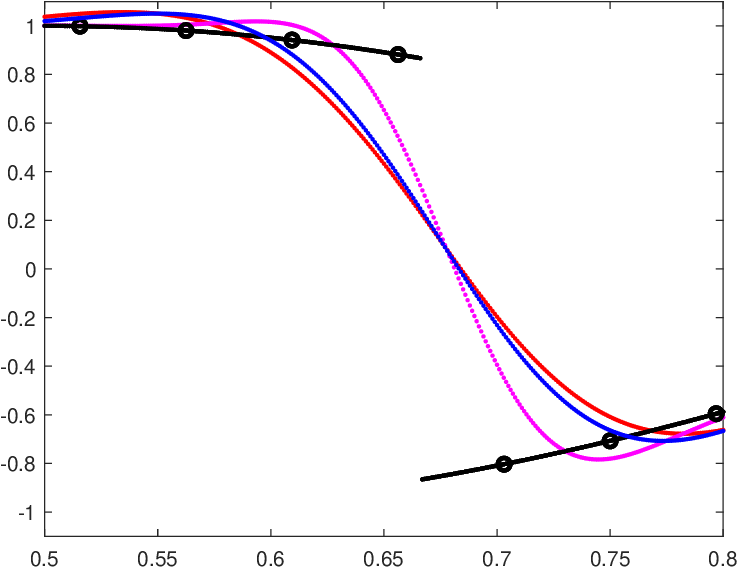}& \includegraphics[width=8cm, height=4.15cm]{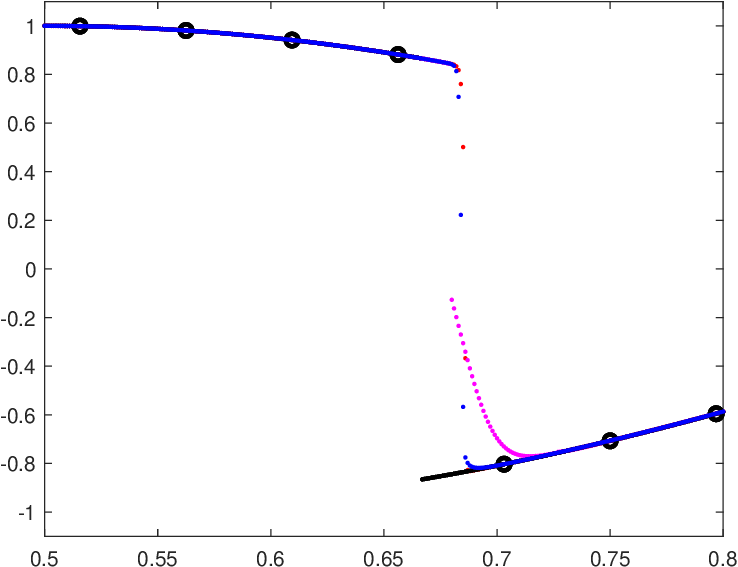}	\\
\multicolumn{2}{c}{$p=3$} \\
\includegraphics[width=8cm, height=4.15cm]{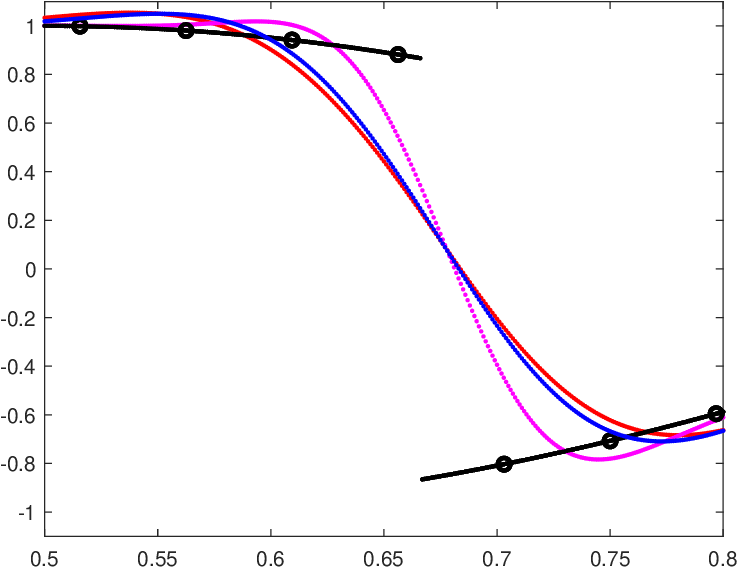}& \includegraphics[width=8cm, height=4.15cm]{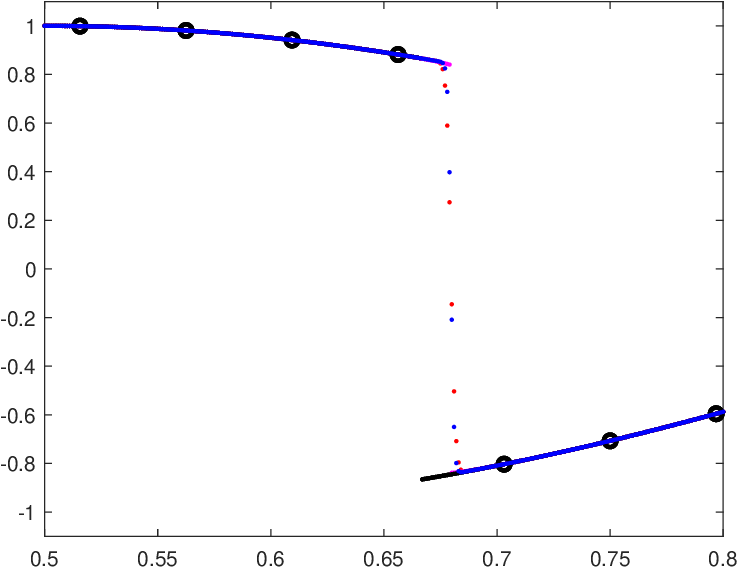}	\\

\end{tabular}
\end{center}
			\caption{Zoom around the discontinuity of the approximation to the function $g$ (black), Eq. \eqref{funciong}, using the linear and non-linear methods with ${\text{W2}}$ (red), ${\text{W4}}$ (blue) and G (magenta).}
		\label{aproxfunciongw2w4gampliacion}
	\end{figure}

\begin{figure}[hbpt!]
\begin{center}
		\begin{tabular}{cc}
	 Linear & 	Non-linear \\
	\includegraphics[width=8cm, height=4.15cm]{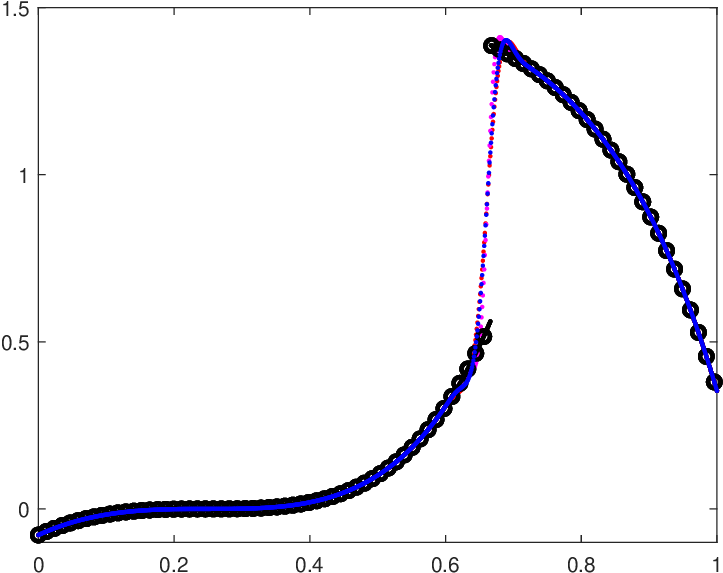} & 	\includegraphics[width=8cm, height=4.15cm]{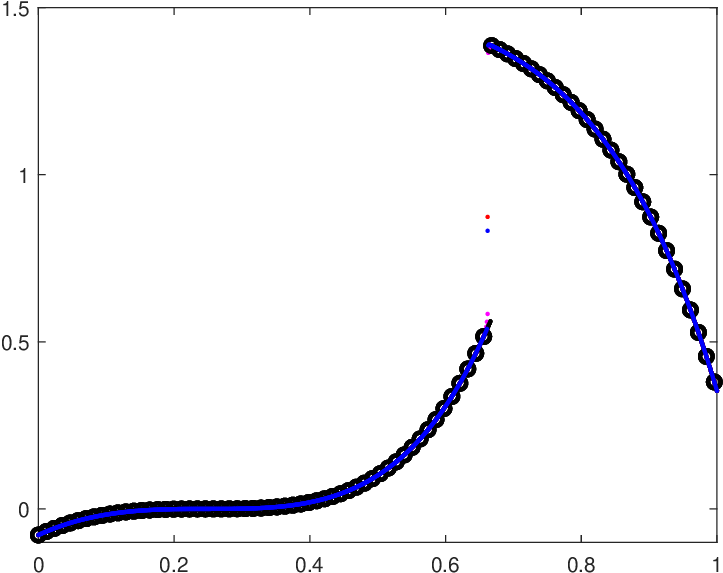}\\
	\includegraphics[width=8cm, height=4.15cm]{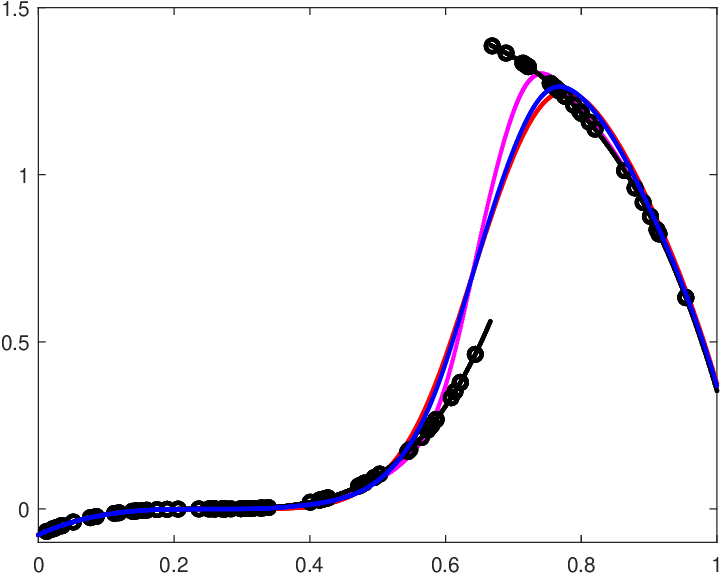} & 	\includegraphics[width=8cm, height=4.15cm]{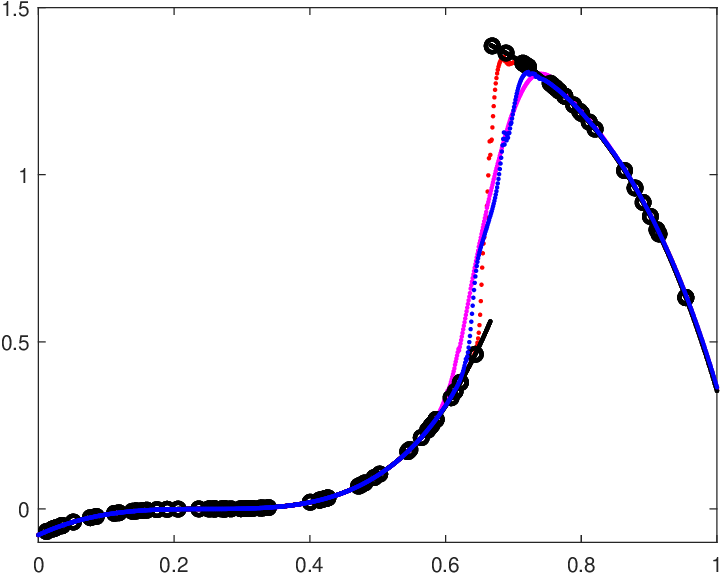}\\
		\end{tabular}
\end{center}
			\caption{Approximation to function $z$ (black), Eq. \eqref{funcionh}, using linear and non-linear methods with ${\text{W2}}$ (red) and ${\text{W4}}$ (blue) and G (magenta); $p=3$. }
		\label{aproxfuncion2w2w4g}
	\end{figure}

In this subsection, we perform two experiments starting with the function:
\begin{equation}\label{funciong}
g(x)=\begin{cases}
\sin(\pi x), & x\leq  2/3, \\
-\sin(\pi x), & x> 2/3, \\
\end{cases}
\end{equation}
discretized in $[-3,3]$ with an uniform grid (rows first and third of Fig. \ref{aproxfunciongw2w4}, Fig. \ref{aproxfunciongw2w4gampliacion}  and first row of Fig. \ref{aproxfuncion2w2w4g}) and with a non-uniform grid $\{x_i\}_{i=0}^{512}$ with $x_i\sim U[-3,3]$, being $U[-3,3]$ a uniform distribution in $[-3,3]$ (second row of Figs. \ref{aproxfunciongw2w4}  and  \ref{aproxfuncion2w2w4g}). The rest of the parameters are identical to those in the previous section. In the left column of Figures \ref{aproxfunciongw2w4}, \ref{aproxfunciongw2w4gampliacion}, and \ref{aproxfuncion2w2w4g}, we present the approximations obtained using the linear methods, while the right column shows the results obtained using the non-linear methods. In all the figures, the original function is plotted with a solid black line, the data points are shown as hollow black circles, the approximation using the function W2 is in red, W4 is in blue, and G is in magenta.%In all the figures, we plot the original function in solid black line, the data are represented through hollow black circles, the approximation using the function W2 in red, W4 in blue and finally, G in magenta. 

%In the column to the left in Figures \ref{aproxfunciongw2w4}, \ref{aproxfunciongw2w4gampliacion} and \ref{aproxfuncion2w2w4g},  we present the approximations obtained when the linear methods are employed, and on the right, the results obtained by the non-linear methods are shown.

First of all, we can observe that
the results obtained are similar independently of the value $p$ chosen, $p=2$ (two first rows), $p=3$ (two last rows), and the weight function $w$. Some non-desirable effects and smearing of the discontinuities appear when linear methods are used, that are avoided or reduced (respectively) when the non-linear
algorithms are employed. When we focus our attention in the results obtained by the non-linear methods, we notice that
we achieve an accurate approximation even in the intervals next to the one containing the discontinuity. This observation is more clearly illustrated in Figure \ref{aproxfunciongw2w4gampliacion}, where we present a zoom of the results around the discontinuity.

Finally, we replicate the same experiment with the function introduced in \cite{ABM}:
\begin{equation}\label{funcionh}
z(x)=\left\{\begin{array}{ll}
5(x-0.25)^3e^{x^2}, & x\leq 2/3,\\
1.5-(x-0.25)^3e^{x^2},&2/3< x.
\end{array}
\right.
\end{equation}
The results for this function are presented in Figure \ref{aproxfuncion2w2w4g}. We can see that we can reach the same conclusions as for the previous experiment. In this case, as the results are analogous for $p=2$ and $p=3$, we only illustrate this last case. In Figure \ref{aproxfuncion2w2w4g}, we present the results for the linear (first column) and the non-linear (second column) methods. The first row presents the results using the uniform grid, and the second row using the non-uniform grid. As before, if we observe Figure \ref{aproxfuncion2w2w4g}, we can see that the new method has a good behavior in all the domain. Near the discontinuity, the oscillations and smearing observed in the results of the linear method have been removed and reduced, respectively, in the non-linear approximation.

\section{Conclusions and future work}\label{conc}

This paper is devoted to adapt the well-known partition of unity method, typically used with radial basis functions, to MLS introducing a non-linear procedure when the data present a discontinuity. The way to insert data dependency in the algorithm is to employ the WENO technique. This method relies in the computation of smoothness indicators to select the data that is free of discontinuities. We have presented the construction of the new algorithm, including the design of smoothness indicators appropriate for our setting. We have also presented some theoretical results regarding the properties of the method. In particular, about the smoothness, polynomial reproduction and order of accuracy. The numerical experiments presented validate the theoretical findings and demonstrate the practical effectiveness of the new method. As future work, it is  important to mention that the parameters $\gamma_k$, $\varepsilon$ (in the case of Gaussian function) and the centers play a key role in getting
an accurate approximation. It is our plan to analyze the  impact of these variables as well as the choice of the function $w$. Additionally, the techniques presented could be adapted for image processing applications in a multiresolution context \cite{harten3}. This would only require modifying the operators $L_i(f)$ and $L(f)$ mentioned in the introduction. Finally, we aim to study the generalization of this non-linear method to multiple dimensions.

\section*{Acknowledgments} 
We would like to thank Dr. David Levin for reviewing various drafts of this paper and for his valuable comments and suggestions during our research. His insights have been invaluable to us.

%We would like to thank David Levin for reading various drafts of this paper, his comments have been very valuable to us.

\section*{Declaration of generative AI and AI-assisted technologies in the writing process}

During the preparation of this work the authors used Microsoft/Copilot to correct orthographical and grammatical errors in the text. After using this tool/service, the authors reviewed and edited the content as needed and take full responsibility for the content of the publication.

%	\begin{figure}[H]
%\begin{center}
%		\begin{tabular}{cc}
%	 Linear & 	Non-linear \\
%\multicolumn{2}{c}{$p=2$} \\
%	\includegraphics[width=8cm, height=5.15cm]{ejemploglinealp2g.eps} & 	\includegraphics[width=8cm, height=5.15cm]{ejemplognolinealp2g.eps}\\
%	\includegraphics[width=8cm, height=5.15cm]{ejemploglinealp2gnoequi.eps} & 	\includegraphics[width=8cm, height=5.15cm]{ejemplognolinealp2gnoequi.eps}\\
%\multicolumn{2}{c}{$p=3$}\\
%	\includegraphics[width=8cm, height=5.15cm]{ejemploglinealp3g.eps} & 	\includegraphics[width=8cm, height=5.15cm]{ejemplognolinealp3g.eps}\\
%	\includegraphics[width=8cm, height=5.15cm]{ejemploglinealp3gnoequi.eps} & 	\includegraphics[width=8cm, height=5.15cm]{ejemplognolinealp3gnoequi.eps}\\
%		\end{tabular}
%\end{center}
%			\caption{Approximation to function $g$, Eq. \eqref{funciong}, using linear and non-linear methods with the gaussian function as $w$.}
%		\label{aproxfunciongg}
%	\end{figure}

\end{document}